# Remarks on a Class of Solutions to the Minimal Surface System

Mu-Tao Wang

## 1. Introduction

The study of the minimal surface equation

(1.1) $$div(\frac{\nabla f}{\sqrt{1+|\nabla f|^2}}) = 0 \, , \, f: D \subset \mathbb{R}^n \to \mathbb{R}$$

is responsible for the progress of nonlinear elliptic PDE theory in the last century. Indeed, most early works on nonlinear elliptic problems focused on this particular equation. There are many beautiful existence, uniqueness, regularity theorems for the minimal surface equation, see for example [13].

The graph of a solution to (1.1) is naturally a minimal hypersurface in $\mathbb{R}^{n+1}$. In general, we can consider a vector-valued function whose graph is a minimal submanifolds of the Euclidean space; the function then satisfies a nonlinear elliptic system. Indeed, a $C^2$ vector-valued function $f = (f^1, \cdots f^m) : D \to \mathbb{R}^m$ is said to be a solution to the minimal surface system (see Osserman [16] or Lawson-Osserman [13]) if

(1.2) $$\sum_{i,j=1}^n \frac{\partial}{\partial x^i}(\sqrt{g} g^{ij} \frac{\partial f^\alpha}{\partial x^j}) = 0 \text{ for each } \alpha = 1 \cdots m$$

where $g^{ij} = (g_{ij})^{-1}$, $g_{ij} = \delta_{ij} + \sum_{\beta=1}^m \frac{\partial f^\beta}{\partial x^i} \frac{\partial f^\beta}{\partial x^j}$, and $g = \det g_{ij}$.

The graph of $f$, consisting of all points $(x^1, \cdots, x^n, f^1(\mathbf{x}), \cdots, f^m(\mathbf{x}))$ with $\mathbf{x} = (x^1, \cdots, x^n) \in D$, is then a minimal submanifold of $\mathbb{R}^{n+m}$ of dimension $n$ and codimension $m$.

The minimal surface system was first studied in Osserman [16], [17] and Lawson-Osserman [13]. In contrast to the codimension one case, Lawson and Osserman [13] discovered remarkable counterexamples to the existence, uniqueness and regularity of solutions to the minimal surface system in higher codimension. It is thus interesting to identify natural conditions under which theorems for the minimal surface equation can be generalized. In this article, we shall discuss a special class of solutions to the minimal surface system. They are vector-valued

2000 *Mathematics Subject Classification.* Primary 53A10, 35J50, 53A07, 49Q05, 53C38.
The author was supported in part by NSF Grant DMS #0104163.

functions that "decrease area" and are natural generalization of scalar functions. After defining area-decreasing maps in the next section, we show several classical results for the minimal surface equation can be generalized within this category.

## 2. Area-decreasing maps

Let $L : \mathbb{R}^n \to \mathbb{R}^m$ be a linear transformation, recall the norm of $L$, $|L|$ is defined by

$$|L| = \sup_{|v|=1} |L(v)|.$$

$L$ induces a liner transformation, $\wedge^2 L$, from the wedge product $\mathbb{R}^n \wedge \mathbb{R}^n = \wedge^2 \mathbb{R}^n$ to $\wedge^2 \mathbb{R}^m$ by

$$(\wedge^2 L)(v \wedge w) = L(v) \wedge L(w).$$

To write $\wedge^2 L$ more explicitly, we pick a basis $\{v_i\}_{i=1\cdots n}$ for $\mathbb{R}^n$ and a basis $\{u_\alpha\}_{\alpha=1\cdots m}$ for $\mathbb{R}^m$, then $\{v_i \wedge v_j\}$, with $(i,j)$ ranges all double indexes $i < j$, forms a basis for $\wedge^2 \mathbb{R}^n$. Likewise $\{u_\alpha \wedge u_\beta\}_{\alpha<\beta}$ forms a basis for $\mathbb{R}^m$. Suppose $L(v_i) = L_{i\alpha} u_\alpha$, in these bases $\wedge^2 L$ is represented by

$$(\wedge^2 L)(v_i \wedge v_j) = L(v_i) \wedge L(v_j) = \sum_{\alpha,\beta} L_{i\alpha} L_{j\beta} u_\alpha \wedge u_\beta = \sum_{\alpha<\beta} (L_{i\alpha} L_{j\beta} - L_{i\alpha} L_{j\beta}) u_\alpha \wedge u_\beta.$$

With this we define

$$|\wedge^2 L| = \sup_{|v \wedge w|=1} |(\wedge^2 L)(v \wedge w)|.$$

In particular, $|\wedge^2 L| = 0$ if $L$ is of rank one.

The norms $|L|$ and $|\wedge^2 L|$ can be expressed by the singular values of $L$, or the eigenvalues of $(L)^T L$. If we denote the singular values by $\lambda_i$, then there exist orthonormal bases $\{v_i\}$ and $\{u_\alpha\}$ such that

$$L(v_i) = \lambda_i u_i$$

if $i$ is less than or equal to the rank of $L$. It is now easy to see

$$|L| = \sup_i \lambda_i$$

and

$$|\wedge^2 L| = \sup_{i<j} \lambda_i \lambda_j.$$

$|\wedge^2 L|$ can also be interpreted as the maximum of the Jacobian of $L$ when restricted to any two-dimensional subspace of $\mathbb{R}^n$.

For a vector-valued function $f : D \to \mathbb{R}^m$, the differential of $f$, $df(x)$ at each $x \in D$ is a linear transformation from $\mathbb{R}^n$ to $\mathbb{R}^m$ represented by the matrix $\frac{\partial f^\alpha}{\partial x^i}$.

DEFINITION 2.1. A Lipschitz map $f : D \to \mathbb{R}^m$ is said to be *area-decreasing* if $|\wedge^2 df|(x) < 1$ for almost every $x \in D$.

In particular, any scalar function, i.e. $m = 1$, is an area-decreasing map.

## 3. Recent results

In this section, we discuss some recent results on area-decreasing minimal maps.

The well-known Bernstein theorem asserts any complete minimal surface that can be written as the graph of an entire function on $\mathbb{R}^2$ must be a plane. This result has been generalized to $\mathbb{R}^n$, for $n \leq 7$ and general dimension under various growth condition, see [3] [5] [20] and the reference therein for the codimension one case. In particular, the following classical result ( see [2], [4], [15]) is an important case of the Bernstein theorem.

THEOREM 3.1. *Let $f : \mathbb{R}^n \to \mathbb{R}$ be a entire solution to (1.1). If $f$ has bounded gradient then $f$ is a linear function.*

For higher codimension Bernstein type problems, there are general results of Simons [21], Reilly [18], Barbosa [1], Ficher-Colbrie [7], Hildebrandt-Jost - Widman,[9] and Jost-Xin[11]. Because of the example of nonparametric minimal cone discovered by Lawson and Osserman [13], all these results seek for optimal conditions under which Bernstein type theorems hold. In [27], we prove the following theorem that improves all known results in higher codimensions. This theorem also generalizes Theorem 3.1 as any scalar function is area-decreasing.

THEOREM 3.2. *Let $f : \mathbb{R}^n \to \mathbb{R}^m$ be a entire solution to (1.2). If $f$ has bounded gradient and $|\wedge^2 df| < 1 - \epsilon$ for some $\epsilon > 0$ then $f$ is a linear map.*

We remark that this is still not the sharpest condition in view of Lawson-Osserman's example.

If we require the graph of $f$ to be a Lagrangian submanifold, (1.2) reduces to a fully nonlinear scalar equation for the potential of $f$. Bernstein type results for such minimal submanifolds have been obtained by many authors, we refer to the recent articles of Yuan[32] and Tsui-Wang [24] and the reference therein.

A corollary of the Bernstein-type result is the following regularity theorem proved in [30].

COROLLARY 3.3. *If $f : D \to \mathbb{R}^m$ is a Lipschitz solution to (1.2) and $|\wedge^2 df| < 1$, then $f$ is smooth.*

The proof of the Bernstein type theorem is based on calculating the Laplacian of the following quantity:

$$*\Omega = \frac{1}{\sqrt{\det(\delta_{ij} + \sum_{\beta=1}^{m} \frac{\partial f^\beta}{\partial x^i} \frac{\partial f^\beta}{\partial x^j})}}.$$

We notice $\sqrt{\det(\delta_{ij} + \sum_{\beta=1}^{m} \frac{\partial f^\beta}{\partial x^i} \frac{\partial f^\beta}{\partial x^j})}$ is indeed the volume element of the graph of $f$.

When $m = 1$, i.e. the codimension one case, $*\Omega = \frac{1}{\sqrt{1+|\nabla f|^2}}$ is the angle made by the normal vector of the graph of $f$ and the $x^{n+1}$ axis. This quantity satisfies a very nice identity in the codimension one case.

$$\Delta \frac{1}{\sqrt{1+|\nabla f|^2}} + \frac{1}{\sqrt{1+|\nabla f|^2}} |A|^2 = 0$$

where $\Delta$ is the Laplace operator of the induced metric on the graph of $f$ and $A$ is the second fundamental form.

This is an important identity that has been used in the codimension-one Bernstein problem, see for example [5] in a slightly different form. In the higher codimension case, we derive in [27] [29]( see also [26]) the following equation for $\ln *\Omega$:

$$(3.1) \qquad \Delta(\ln *\Omega) = -\sum_{\alpha,l,k} h_{\alpha lk}^2 - \sum_{i,j,k} \lambda_i \lambda_j h_{n+i,jk} h_{n+j,ik}$$

where $\lambda_i$ are singular values of $df$ and $h_{\alpha lk}$ are coefficients of the second fundamental form in special bases adapted to the singular valued decomposition of $df$. In particular, $\sum_{\alpha,l,k} h_{\alpha lk}^2 = |A|^2$ is the square norm of the full second fundamental form.

The right hand side can be rewritten as

$$-\sum_{\alpha,l,k} h_{\alpha lk}^2 - \sum_{i,k} \lambda_i^2 h_{n+i,ik}^2 - 2\sum_{i<j,k} \lambda_i \lambda_j h_{n+i,jk} h_{n+j,ik}.$$

It is clear now if $\lambda_i \lambda_j < 1$ for $i \ne j$, or $f$ is area-decreasing, the righthand side is always non-positive by completing square. The proof of Theorem 3 is based on blow-up analysis of this inequality.

Next we discuss the interior gradient bound for solutions to the minimal surface system (1.2). For solutions to (1.1), this was discovered, in the case of two variables, by Finn [6] and in the general case by Bombieri, Di Giorgi and Miranda [2]. The a priori bound is a key step in the existence and regularity of minimal surface theory. The estimate has been generalized to other curvature equations and proved by different methods. We refer to the note at the end of chapter 16 of [8] for literature in these directions and Theorem 16.5 of [8] for the precise statement in the codimension-one case. In [29], we generalize the interior gradient estimate to higher codimension:

THEOREM 3.4. *Let $D$ be a domain in $\mathbb{R}^n$ and $f : D \to \mathbb{R}^m$ a $C^2$ solution to equation (1.2) such that $|\wedge^2 df| < 1$. If each $f^\alpha$ is non-negative, then for any point $x_0 \in D$, we have estimate*

$$|df(x_0)| \le C_1 \exp\{C_2 |f(x_0)|/d\}$$

*where $d = dist(x_0, \partial D)$ and $C_1, C_2$ are constants depending on $n$.*

The proof of Theorem 5 again uses (3.1). Indeed, if $\lambda_i \lambda_j \le 1$ for $i \ne j$ then we derive

$$\Delta \ln *\Omega \le -\frac{1}{n} |\nabla \ln *\Omega|^2$$

where $\nabla$ denotes the gradient of the induced metric on the graph of $f$. In [29], we provide two arguments to derive the gradient estimate, the first one follows the approach developed by Michael-Simon [14] and Trudinger [22][23] as was presented in [8] and the second one generalizes Korevvar's [12] using the maximum principle (see also [31]).

The underlying principle in deriving (3.1) is discussed in [28]. It turns out $|\lambda_i \lambda_j| < 1$ for $i \ne j$ defines a region on the Grassmannian of $n$-planes in $\mathbb{R}^{n+m}$ on which $-\ln \Omega$ is a convex function. That $\ln *\Omega$ is superharmonic follows directly from Ruh-Vilms' [19] theorem which states the Gauss map of a minimal submanifold is always a harmonic map.

REMARKS ON A CLASS OF SOLUTIONS TO THE MINIMAL SURFACE SYSTEM    5To get a feeling of what the region $|\lambda_i \lambda_j| < 1$ looks like on the Grassmannian, we demonstrate the case $n = m = 2$ in the following . In this case, $G(2,2) = S^2(\frac{1}{\sqrt{2}}) \times S^2(\frac{1}{\sqrt{2}})$, where $S^2(\frac{1}{\sqrt{2}})$ is a two-sphere of radius $\frac{1}{\sqrt{2}}$. Write $\mathbb{R}^4 = \mathbb{R}^2 \oplus \mathbb{R}^2$ and denote the coordinates on the first summand $\mathbb{R}^2$ by $x^1, x^2$ and the coordinates on second summand $\mathbb{R}^2$ by $y^1, y^2$. Then the forms $\omega_1$ and $\omega_2$

$$\omega_1 = \frac{1}{\sqrt{2}}(dx^1 \wedge dx^2 - dy^1 \wedge dy^2)$$

$$\omega_2 = \frac{1}{\sqrt{2}}(dx^1 \wedge dx^2 + dy^1 \wedge dy^2)$$

are viewed as functions on $G(2,2)$, see for example section 3 in [**28**]. We may identify $\omega_1$ and $\omega_2$ with the height functions on the first and the second factor in $S^2(\frac{1}{\sqrt{2}})$, respectively. For a linear transformation $L$, the condition $|\lambda_1 \lambda_2| < 1$ corresponds to $|dy^1 \wedge dy^2| < dx^1 \wedge dx^2$ or $\omega_1 > 0$ and $\omega_2 > 0$ on the graph of $L$. Therefore, the area-decreasing condition on $f$ is equivalent to the image of the Gauss map of the graph of $f$ lies in the product of two hemispheres of $S^2(\frac{1}{\sqrt{2}}) \times S^2(\frac{1}{\sqrt{2}})$. Recall in the codimension-one case, a hypersurface in $\mathbb{R}^{n+1}$ is the graph of a scalar function $f$ if and only the Gauss map lies in a hemisphere of $S^n$. This also indicates why the area-decreasing condition is a natural generalization of being codimension-one. We conjecture most classical theorems for solutions to the minimal surface equation should hold for this class of maps. In particular,

CONJECTURE 3.5. Let $D$ be a $C^2$ convex domain in $\mathbb{R}^n$. If $\phi : \overline{D} \to \mathbb{R}^m$ is an area-decreasing map, then the Dirichlet problem of the minimal surface system (1.2) for $\phi|_{\partial D}$ is solvable in smooth maps. The solution is volume-minimizing and is unique among all area decreasing maps with the same boundary condition.

The Dirichlet problem for the minimal surface equation (1.1) was solved by Jenkin-Serrin [**10**] for mean convex domains. The solution of Conjecture 6 will be a natural generalization of Jenkin-Serrin's theorem to systems. In [**30**], we use the mean curvature flow to prove the existence part of conjecture under a stronger condition on $\phi$. It is shown in [**28**] that the set of area-decreasing linear transformations is a convex subset of the Grassmannian; we believe this will be useful in attacking Conjecture 6.

## References

[1] Barbosa, Jo ao Lucas Marqus, *An extrinsic rigidity theorem for minimal immersions from $S^2$ into $S^n$.*, J. Differential Geom. 14 (1979), no. 3, 355–368 (1980).
[2] E. Bombieri, E. De Giorgi and M. Miranda, *Una maggiorazione a priori relativa alle ipersuperfici minimali non parametriche.* (Italian) Arch. Rational Mech. Anal. 32 (1969) 255–267.
[3] L. Caffarelli, L. Nirenberg, and J. Spruck, *On a form of Bernstein's theorem. Analyse mathmatique et applications.* 55–66, Gauthier-Villars, Montrouge, (1988).
[4] E. De Giorgi, *Sulla differenziabilit e l'analiticit delle estremali degli integrali multipli regolari.* (Italian) Mem. Accad. Sci. Torino. Cl. Sci. Fis. Mat. Nat. (3) 3 (1957) 25–43.
[5] K. Ecker and G. Huisken, *A Bernstein result for minimal graphs of controlled growth.*, J. Differential Geom. 31 (1990), no. 2, 397–400.
[6] R. Finn, *On equations of minimal surface type.* Ann. of Math. (2) 60, (1954). 397–416.
[7] D. Fischer-Colbrie, *Some rigidity theorems for minimal submanifolds of the sphere.*, Acta Math. 145 (1980), no. 1-2, 29–46.

Department of Mathematics, Columbia University, New York, NY 10027
*E-mail address*: mtwang@math.columbia.edu